\documentclass[preprint,12pt]{elsarticle}

    \makeatletter
    \def\ps@pprintTitle{%
      \let\@oddhead\@empty
      \let\@evenhead\@empty
      \def\@oddfoot{\reset@font\hfil\thepage\hfil}
      \let\@evenfoot\@oddfoot
    }
    \makeatother

\usepackage{amsmath}
\usepackage{array}
\usepackage{url}

\usepackage{graphicx}
\usepackage{epstopdf} 

\usepackage{amssymb}

\begin{document}

\begin{frontmatter}

\title{Exhaustive generation of `Mrs Perkins's quilt' square dissections for low orders}

\author{Ed Wynn}

\address{175 Edmund Road, Sheffield S2 4EG, U.K.}
\ead{ed.wynn@zoho.com}

\begin{abstract}
Dissections of a square into smaller squares, with the smaller squares having relatively prime sizes, are known as Mrs Perkins's quilts.  A representation of these dissections using graphs is presented.  The edges are directed and  coloured North-South or West-East, and the graph corresponds naturally to the dissection.  This representation allowed the exhaustive generation of all dissections up to order 18, using the plantri software.  The results were cross-checked by generating all dissections of small sizes using a direct approach.  The results confirm, extend and introduce several  integer sequences.

\end{abstract}


\end{frontmatter}


\section{Introduction}
\label{introduction}
We are interested in dissecting squares into smaller squares; the \emph{order} of a dissection is the number of smaller squares, which are here called \emph{subsquares}.  We do not require that the dissection be \emph{perfect} (that is, that the subsquares are all of distinct sizes), but we do require that it is \emph{prime} (that is, that the greatest  divisor of their sizes is 1).  Such dissections are known as \emph{Mrs Perkins's quilts}, after the name of a puzzle \cite[Problem 173]{Dudeney1917}.

We define the function $f(n)$ to be the least possible order of a prime dissection of a square of side $n$.
It is known that $\log_{2}n < f(n) < 6 \log_{2}n$ for $n>1$; the lower bound is due to Conway \cite{Conway1964} and the upper bound is due to Trustrum \cite{Trustrum1965}.  Determination of $f(n)$ is problem C3 in \cite{CroftFalconerGuy1991}, where values were supplied `with some slightly increasing lack of confidence' for $n\leq100$.  The values came from \cite{Conway1964}, where it was admitted: `Strictly speaking, these are only upper bounds for $f(n)$. Readers may care to accept the implied challenge.'  The lack of confidence was justified: lower orders were subsequently found for several values of $n$, starting with 51.  The best values known to date are listed in the On-Line Encyclopedia of Integer Sequences \cite[{S}equence {A}005670]{OEIS}.  However, some of the uncertainty persists: in the OEIS entry, it is stated that `$n\leq15$ (and possibly 16) proved minimal by J. H. Conway'.\footnote{The OEIS entry cites an entry on the math-fun mailing list, `Re: [math-fun] {M}rs. {P}erkins {Q}uilt - {O}rders 89, 90 improved over {UPIG}.', 10 October 2003.  I am grateful to Stuart Anderson for providing a copy of Prof. Conway's message: `Let me say that there should have been no lack of confidence up to 15, because I proved the above numbers best possible up to there long ago (I think in fact up to 16, but am not quite sure of that)'.}

The current article reports an exhaustive generation of prime dissections of squares with orders up to 18.  This confirms the currently conjectured values of $f(n)$ for $n\leq120$.  Several other sequences are extended.

The current approach is to define a way in which dissections can be represented by graphs, as discussed in the next section.  The representation used here can be compared to the use of graphs in similar areas -- most notably, the analogy with electrical networks.  One of the inventors of that analogy (in 1940) has more recently written an entertaining account of that work \cite{Tutte1958}.  The electrical-network analogy (here abbreviated to ENA) was successfully applied to discovery of perfect dissections \cite{Duijvestijn1978}.  The graph in that approach is a representation of a scan in either the horizontal or the vertical direction.  The vertices correspond to boundaries between subsquares, and an edge connecting two vertices corresponds to a subsquare that spans those boundaries.  The network can be efficiently solved to give sizes of subsquares.  Planar imbeddings of these graphs give possible layouts of dissections, but these are not guaranteed to fit together in the other direction, or (if they do fit together) to form squares.  For perfect dissections, the graphs have minimum degree 3, but this simplification would not be appropriate if they were used to investigate imperfect dissections.  Reviews in this area, concentrating on perfect dissections, are by Federico \cite{Federico1979} and Anderson \cite{Anderson2013}.

The alternative approach described here also uses planar imbeddings of graphs to represent the layout of dissections.  One difference in the graphs is that, as described in Section~\ref{graph}, vertices here represent subsquares, rather than boundaries between subsquares.  Another difference is that here a graph includes both horizontal and vertical connections, so a set of sizes consistent with a graph will fit together to a full square.  Considerably more computation is required for each graph, including a search for ways to direct and colour it, but the requirements of triangulation and minimum degree 4 reduce the number of candidate graphs.  Despite the differences between the graphs, there are considerable similarities between the equations generated; these are discussed in Section~\ref{examplecalcs}.

The graph representation used here is conceptually straightforward: vertices represent subsquares, and edges represent connections.  Therefore, similar graphs have been used by several authors, and applied to architectural floor plans \cite{Mitchell1976}, for example.  Several variations are discussed by Felsner \cite{Felsner2012}.  In particular, Fusy \cite{Fusy2009} investigates irreducible triangulations of the 4-gon with transversal structures (which are defined in very similar terms to the direction and colouring of the edges, described in the next section).  Fusy proved several properties of these graphs, and applied the results to straight-line drawings of planar graphs.  However, it has been shown \cite{Felsner2012} that the results are also relevant to square or rectangle dissections.

In contrast to the methods developed to solve specific problems, it is interesting to note an alternative approach \cite{Kurz2012}, where the problem is stated in terms of integer linear programming.  This can then be solved using general solvers.  This approach may reduce the need for ingenuity and thought; the extra expense in computer time, if any, has not yet been fully investigated.

\section{Representation of dissections by two-coloured, directed graphs}
\label{graph}

\begin{figure}
\centering
\includegraphics[width=8cm]{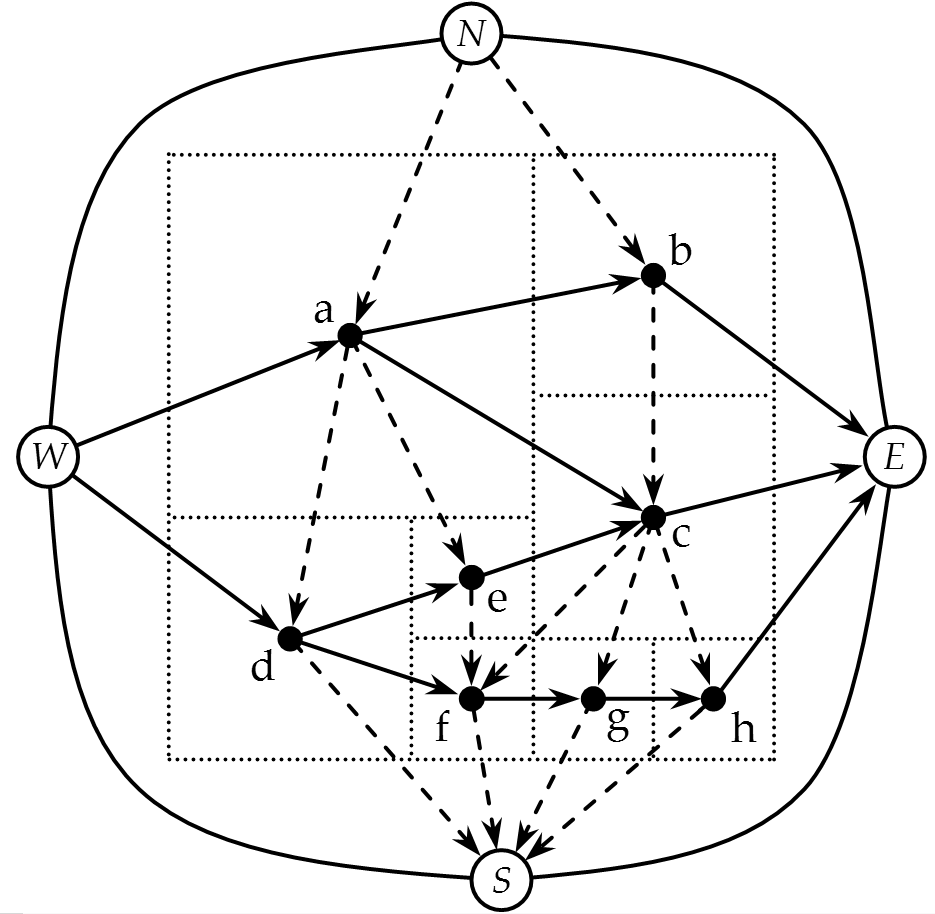}
\caption{An example of a dissection (dotted lines) overlaid with the corresponding graph.
The four circles are cardinal vertices, referred to as $V_N$ etc in the text.  North-South edges are dashed.  The lowercase labels for subsquares are used in Section~\ref{examplecalcs}.}
\label{fig:graph}
\end{figure}

The basis for the generation of the dissections is their representation by graphs.  An example is shown in Figure~\ref{fig:graph}.  Let $D$ be a dissection, and $G$ be a graph that represents $D$ in the following way:
\begin{itemize}
\item Each subsquare in $D$ corresponds to a vertex in $G$.
\item In addition, there are four \emph{cardinal} vertices, so called because they can be labelled North, East, South and West in a suitable cyclic order.  These will be referred to as $V_N$, $V_E$, $V_S$ and $V_W$ respectively, and lower-case $v$ will always refer to non-cardinal vertices.
\item Where two subsquares in $D$ are contiguous, the corresponding vertices are connected.  Each edge is directed and has one of two colours, North-South and West-East, depending on the relative placement of the subsquares.
\item Subsquares on the four extremes of the dissected square are connected to the corresponding cardinal vertices, again with edges directed and coloured in the natural way.
\item Where four subsquares meet at a \emph{cross} (that is, where four subsquares have a point in common), exactly one of the two diagonal pairs of subsquares is connected.
  \begin{itemize} \item This feature ensures that $G$ is planar and triangulated.  However, it implies that dissections with crosses can be represented by more than one graph: for each cross, there is a choice of which pair to connect, and which  colour the edge should have.  There is one cross in Figure~\ref{fig:graph}, involving three unit subsquares and one $2\times2$ subsquare.  The latter has been connected to its diagonal partner with a North-South edge, but a West-East edge with the opposite direction is an alternative, and so are two directions in the other diagonal connection. \end{itemize}
\item Each cardinal vertex is connected to the two neighbouring cardinals in the cycle.  The direction and colour of these edges is arbitrary; to be consistent with the properties discussed below, they can be North-South edges, directed away from North or towards South.
  \begin{itemize} \item  These edges ensure that all vertices have degree at least 4 (as discussed below). \end{itemize}
\end{itemize}

Properties of $G$ can be deduced:
\begin{itemize}
\item If $D$ has order $N$, then $G$ has order $N+4$.
\item $G$ is planar, and the imbedding follows the same layout as $D$.
\item Every face of $G$ is a triangle, except the infinite outer face; this face has four edges which connect the cardinal vertices.  The faces are triangles because in the dissection there is no gap between two adjacent subsquares that both touch a third subsquare.  No triangle contains other vertices.
\item In every triangle, two edges have the opposite colour from the third, and those two edges either both point towards their shared vertex or both point away from it.
\item Every $v$ has the following edges, in clockwise order: at least one North-South edge, pointing to the vertex; at least one West-East edge, pointing away; at least one North-South edge, pointing away; and then at least one West-East edge, pointing to the vertex.  This is because each subsquare in the dissection has at least one neighbour on each face, if we include the cardinal vertices for subsquares on the extremes.  It follows that every non-cardinal vertex has degree at least 4.
\item No $v$ is connected to both $V_N$ and $V_S$, or to both $V_W$ and $V_E$.  Such a subsquare would need to have length equal to that of the dissected square; the search does not include the \emph{trivial prime dissection} of a unit square dissected into a single unit subsquare.  It follows that at least two non-cardinal vertices are connected to each $V$, and that every $V$ has degree at least 4.  These properties would not apply in general if this method were applied to dissections of rectangles into smaller rectangles.
\item Exactly one $v$ is connected to each other pair of cardinal vertices; this is called a \emph{corner} vertex.
\item No edge connects $V_N$ to $V_S$, nor $V_W$ to $V_E$.
\item Edges from $V_N$ to non-cardinal vertices are coloured and directed in the natural way: North-South, pointing away from $V_N$.  These vertices are linked by West-East edges in a chain from the North-West corner to the North-East corner, in anticlockwise order of their connection to $V_N$.  There are no other connections between these vertices.  Similar situations apply for the other cardinal vertices.
\end{itemize}

The current work used the plantri software \cite{plantri} to generate an exhaustive set of candidate graphs; in the terminology of this software, we must consider one member of each isomorphism class of imbedded triangulations of a 4-sided disk.  Plantri is well suited to this task.   The plantri software was adapted to restrict the search to graphs with minimum degree 4; this restriction was not available by default for imbeddings of a disk.  Some candidates could be rejected on the basis of the required properties, discussed above.  For all other candidates, all possible directions and colourings were considered, using an exhaustive backtracking search.  Each new direction and colouring was tested for conformance with the required properties discussed above; unsuitable choices were rejected.  The remaining directed and coloured graphs were analysed to deduce, if possible, lengths of subsquares that would produce a corresponding dissection of the square.  This analysis is described in the next section; each graph considered was shown to correspond to a single prime dissection of a square, or shown not to correspond to any.  Thus, an exhaustive search was conducted for dissections of small orders; results are presented in Section~\ref{results}.

\section{Deducing sizes}
\label{equations}

From now on, we will consider only coloured, directed graphs, in specified imbeddings, that satisfy the properties mentioned in the previous section.  For a dissection of order $N$, the lengths of the subsquares, $L_1,\dots,L_N$, can be regarded as unknowns, and any such graph can be used to deduce equations that these lengths obey.  For example, define a \emph{traverse} to be a path of North-South edges from $V_N$ to $V_S$, or an equivalent West-East path.  The sum of the lengths associated with the non-cardinal vertices in a traverse equals the length $n$ of the dissected square, so each traverse provides an equation.  Kurz \cite{Kurz2012} gives an example of these equations (not using the graph representation).

At this stage, $n$ is also unknown.  Instead of equating sums of lengths to $n$, we may define a subsquare's \emph{normalised length} to be its length divided by $n$. Each sum of these normalised lengths is then equated to 1.  In a successful solution, the resulting normalised lengths will then be rational numbers (since all the equations have integer, indeed unity, coefficients); they must also satisfy $0<L_i<1$.  If all normalised lengths are multiplied by the least common multiple of all their denominators, the result will be a prime dissection with integer lengths.  For some graphs, the equations have a solution with normalised lengths not satisfying $0<L_i<1$; this indicates that the graph does not correspond to a dissection of the square.  Whether a solution uses lengths or normalised lengths, a multiplicative constant can be found to give a prime dissection; therefore, it is not always necessary to distinguish between lengths and normalised lengths.  It has been noted previously \cite{BrooksSmithStoneTutte1992} that solutions to dissections typically provide lengths fixed to within a multiplicative constant.  As discussed in Section~\ref{examplecalcs}, none of the systems of linear equations had more than one solution in the current work.

The number of traverses, and hence the number of equations, will depend on the graph.  However, in the electrical-network analogy (ENA), even the horizontal or vertical constraints are sufficient to define a minimal integer solution.  (This is noted by Kurz~\cite{Kurz2012}.  Acton~\cite{Acton1990} discusses the solution of similar sets of equations.)  The current equations contain both horizontal and vertical constraints.  Section~\ref{examplecalcs} gives an example of the specific equations used and shows that the different approaches are closely related.

Some dissections (specifically, those containing crosses) have several graphs.  All graphs were generated, and duplicate solutions were removed.

Some solutions were generated that satisfied the equations, but did not match the graph that defined the equations: for example, solutions contained pairs of subsquares that were connected in the graph and were correctly aligned in the deduced dissection (for example, in that the lower edge of one had the same vertical coordinate as the top edge of the other) but did not actually touch, even at a cross.  These solutions were discarded; if the dissections were valid, they would also be generated from their correct graphs.

In fact, the search did not use the approach of building up equations from traverses, because the full collection of traverses is not available until the graph has been fully coloured and directed.  Instead, a greater number of unknowns was used, so that each new direction and colouring of an edge would generate a new equation.  We refer to this as the \emph{local equation} approach.  Each new local equation could then be tested for compatibility with the existing set of equations, and impossible graphs could be rejected at a higher level of the search tree.  This approach used $3N$ unknowns in total: the $N$ lengths of the subsquares, and also the $x$ and $y$ coordinates of the subsquares' North-West corners.  Some of the unknowns were known immediately: for example, with $x$ pointing East, $x=0$ for all subsquares on the Western extreme; these subsquares were known from their connection to $V_W$.  Solving for fractional dissections of a unit square (in other words, using normalised sizes), then a subsquare $i$ on the Eastern extreme satisfied $(x_i+L_i)=1$.  During the direction/colouring search, the assignment of a direction and a colour to an edge produced a new equation immediately.  For example, if a West-East edge pointed from subsquare $i$ to subsquare $j$, then $(x_i+L_i)=x_j$.  The additional complication of $3N$ unknowns was found to be justified by the earlier detection of graphs without solutions as dissections.

Equations were combined by forming them into an upper-triangular matrix, for a specified ordering of the unknowns.  This was effectively an incremental process of Gaussian elimination.  Each new equation was considered in terms of its earliest unknown.  If no equation had previously been noted with that unknown as the earliest, then the new equation contributed new information and could not be inconsistent with the previous equations; the new equation would be duly noted.  Alternatively, if a previous equation already had that unknown as its earliest, then the two equations would be combined to eliminate that unknown.  If the combined equation contained no unknowns, then there were two possible conclusions: if the combined equation contained a non-zero constant term, then the new equation was not consistent with the existing equations; alternatively, if the combined equation effectively stated $0=0$, then the new equation was redundant, since it was a linear combination of previous equations.  If the combined equation had non-zero coefficients of unknowns, then the process would be repeated in terms of the new earliest unknown.

Integer arithmetic was used to combine equations.  It was verified that overflow did not occur.  The magnitudes of coefficients often increased when equations were combined; to keep the coefficients within manageable magnitudes, it was necessary to divide every newly-combined equation by the greatest common divisor of its coefficients.  For order $N\leq18$, all terms had magnitudes less than $2^{43}$, and so 64-bit `long long' integers in the C programming language were adequate.

\section{Relationship between different approaches}
\label{examplecalcs}

This section presents the equations that apply to the dissection of Figure~\ref{fig:graph} using three approaches: the electrical-network analogy (ENA) of \cite{Tutte1958}; the new approach using traverses; and the new approach using local equations.  By reference to the example, it is shown that the approaches are closely related.

\begin{figure}
\centering
\includegraphics[width=8cm]{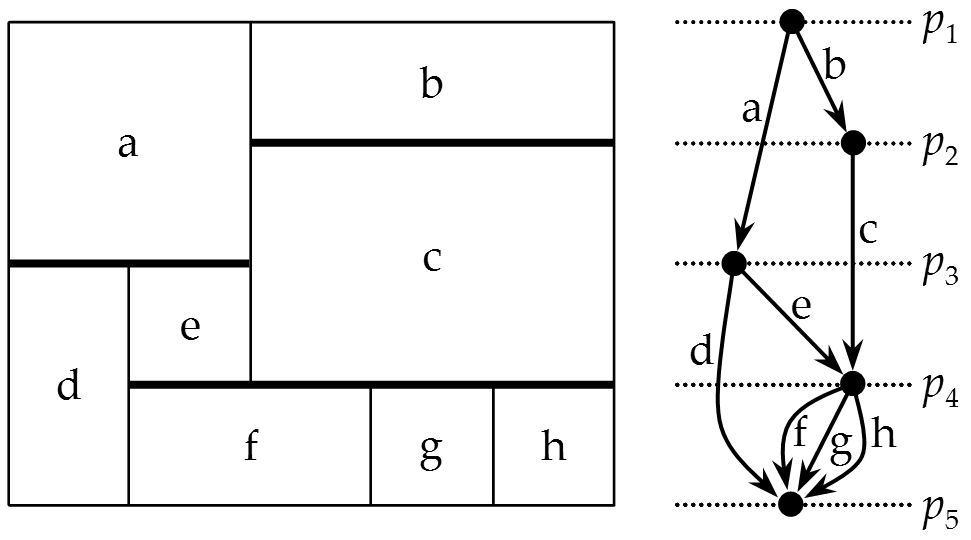}
\caption{The dissection of Figure~\ref{fig:graph}, with the graph from the vertical scan in the electrical-networks analogy.  The dissection is shown without square proportions, reflecting the uncertainty still present before the equations are solved.  As before, the cross between subsquares e, c, f and g has been broken by selecting a vertical connection from c to f.}
\label{fig:electric}
\end{figure}

\subsection{Electrical-network analogy (ENA)}
\label{subsection:electrical}

We start with the ENA and consider only the vertical (North-South) scan, whose graph is shown in Figure~\ref{fig:electric}.  In a formulation by the ENA's original authors \cite{BrooksSmithStoneTutte1992}, a matrix $\mathbf{i}$ is defined as follows:
\begin{equation} \label{definei} i_{hr} = \begin{cases}
    1 & \text{if node $N_h$ is at the top of edge $E_r$} \\
   -1 & \text{if node $N_h$ is at the base of edge $E_r$} \\
    0 & \text{otherwise.}
\end{cases}
\end{equation}
In general let the number of nodes be $m$.  For the graph in Figure~\ref{fig:electric}, with 5 nodes and 8 edges, the matrix is
\begin{equation} \label{definei}  \mathbf{i} = \left( \begin{array}{*{8}r}
 1 & 1 & 0 & 0 & 0 & 0 & 0 & 0 \\
 0 & -1 & 1 & 0 & 0 & 0 & 0 & 0 \\
 -1 & 0 & 0 & 1 & 1 & 0 & 0 & 0 \\
 0 & 0 & -1 & 0 & -1 & 1 & 1 & 1 \\
 0 & 0 & 0 & -1 & 0 & -1 & -1 & -1 
\end{array} \right) .\end{equation}
A matrix $\mathbf{j}$ is defined:
\begin{equation} \label{definej} \mathbf{j} = \mathbf{i} \mathbf{i}^\intercal. \end{equation} 
The equations to be solved are then:
\begin{equation} \label{jmatrix} \sum_{k=1}^{m}{j_{hk} p_k} = 0 \end{equation} 
for all rows of $\mathbf{j}$ except the top and bottom rows: $h\neq1,m$ (using the original formulation where indices start at 1).

The original formulation refers to the unknowns $\mathbf{p}$ as `potentials', but here we will explore the geometric basis for the equations, so we note that $p_k$ is the vertical coordinate of node $N_k$; for example $p_1=y_a=y_b$.  It is observed \cite{BrooksSmithStoneTutte1992} that one of the coordinates, such as $p_1$, can be fixed to be 0 without loss of generality.  Therefore (\ref{jmatrix}) represents $m-2$ equations in the remaining $m-1$ unknowns.  It has been shown \cite{BrooksSmithStoneTutte1992} that these equations are linearly independent, and therefore they fix the remaining coordinates to within a multiplicative constant.

In the example, we will take several steps: form $\mathbf{j}$ from $\mathbf{i}$; remove the top and bottom rows; ignore $p_1$ since it is fixed to be 0; take $p_5$ to the right of the equation; and divide through by $p_5$.  The result is:
\begin{equation} \left( \begin{array}{*{3}r}  2 & 0 & -1 \\ 0 & 3 & -1 \\ -1 & -1 & 5 \end{array} \right) 
\left( \begin{array}{r} \left. p_2 \middle/ p_5 \right. \\\left. p_3 \middle/ p_5 \right. \\\left. p_4 \middle/ p_5 \right.  \end{array} \right)  = 
\left( \begin{array}{r} 0 \\ 1 \\ 3  \end{array} \right)  
\end{equation}
Noting that $p_5$ equals $n$, the length of the square, we have $m-2$ linearly independent equations in $m-2$ normalised coordinates, giving the solutions $2/5$, $3/5$ and $4/5$.  Figure~\ref{fig:graph} confirms that these are correct.

For geometrical insight, it may be preferable to expand $\mathbf{j}$ in terms of $\mathbf{i}$.  Each edge in the ENA's graph represents a subsquare of the dissection, and is represented by a column of $\mathbf{i}$.  This column contains a 1 and a $-1$ for the starting and ending nodes.  Therefore, $\mathbf{i}^\intercal\mathbf{p}$ equals $-\mathbf{L}$, where $\mathbf{L}$ is the vector of subsquares' lengths, possibly normalised.  For example,
\begin{equation} \label{lengthpdiff} -L_a = p_1-p_3 .\end{equation}

Applying the expansion of $\mathbf{j}$ to equation (\ref{jmatrix}) in our example, we obtain:
\begin{equation} \label{electricalequations}
\begin{array}{lr} 
-L_b+L_c & = 0 \\ 
-L_a+L_d+L_e & = 0 \\ 
-L_c-L_e+L_f+L_g+L_h & = 0
\end{array} \end{equation}
The geometric origin of these equations can be found, perhaps surprisingly, from the horizontal traverses in the new approach.  This is discussed in the next section.

\subsection{Traverse equations}
\label{subsection:traverse}

By reference to Figure~\ref{fig:graph}, the traverses in the example can be found.  From North-South traverses, the equations are:
\begin{equation} \label{traversesNS}
\begin{array}{lr} L_a+L_d & = n \\ 
  L_a+L_e+L_f& = n \\ 
  L_b+L_c+L_f & = n \\
  L_b+L_c+L_g & = n \\
  L_b+L_c+L_h & = n 
\end{array} \end{equation}
From West-East traverses, the equations are:
\begin{equation} \label{traversesWE}
\begin{array}{lr} 
  L_a+L_b & = n \\
  L_a+L_c & = n \\
  L_d+L_e+L_c & = n \\
  L_d+L_f+L_g+L_h & = n
\end{array} \end{equation}
Since the right sides of these equations are all $n$, the same left sides can be expressed in terms of normalised lengths and equated to 1.  This is computationally much more convenient, because $n$ is not known at the start.  In fact, $n$ is fixed by the specification of a prime dissection rather than by the geometrical constraints.

The connection with the ENA can be observed by taking the differences of successive equations in (\ref{traversesWE}): the resulting equations are identical to (\ref{electricalequations}).  This is despite the fact that (\ref{electricalequations}) was generated by the \emph{vertical} scan in the ENA, and (\ref{traversesWE}) from \emph{horizontal} traverses.  The reason is that both these equations represent constraints on horizontal distances, in equivalent situations, as will now be discussed.

\begin{figure}
\centering
\includegraphics[height=3.5cm]{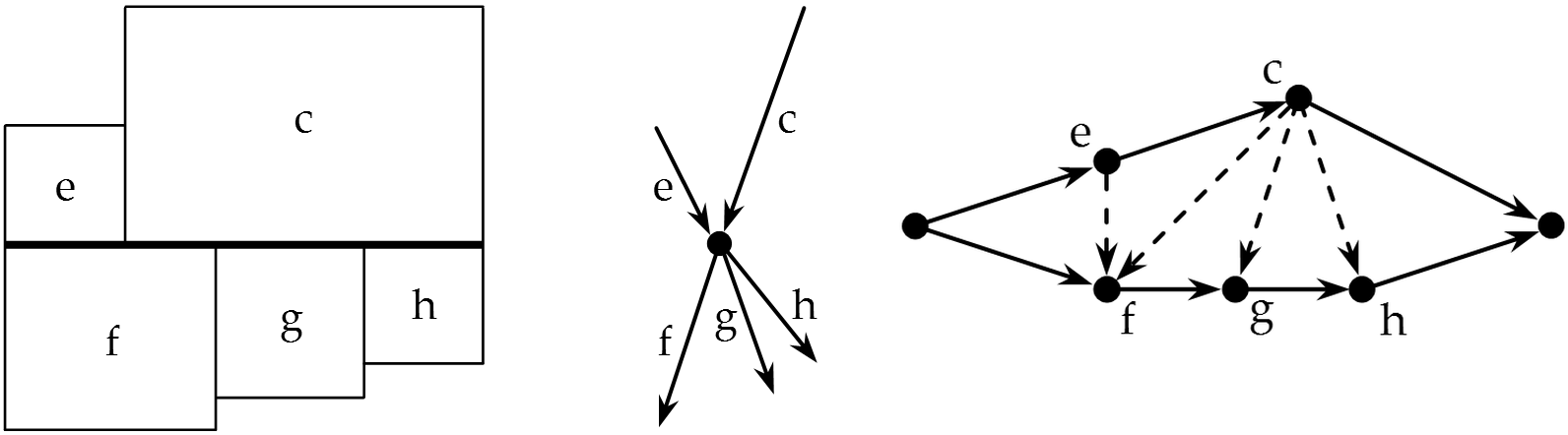}
\caption{Parts of (left to right) the dissection of Figure~\ref{fig:graph}, the vertical ENA graph of Figure~\ref{fig:electric} and the graph of Figure~\ref{fig:graph}.  These contain the same information about sums of subsquares' widths: $L_e+L_c=L_f+L_g+L_h$.  (The dissection and the graph on the right contain the additional information that $L_c\geq L_g+L_h$, because part of subsquare c is vertically above f.)}
\label{fig:level}
\end{figure}

An example of the correspondence is shown in Figure~\ref{fig:level}.  In the ENA's graph, only edges e and c enter the node, and only edges f, g and h leave the node.  It can be deduced that this part of the dissection must be similar to the part shown in the figure: the West edges of e and f are aligned, and so are the East edges of c and h. We deduce that
\begin{equation} \label{deducedlengths} L_c+L_e = L_f+L_g+L_h \end{equation}
which is identical to the last equation in (\ref{electricalequations}).  Thus we have a geometric interpretation of the nodes in the ENA's graph: each node corresponds to a horizontal line in the dissection, extending as far as possible to West and East.  (This interpretation is, of course, not novel, being consistent with the original expositions of this method.)  The equations in (\ref{jmatrix}) follow: $\mathbf{j}$ is expanded using (\ref{definej});  $\mathbf{i}^\intercal\mathbf{p}$ equals $-\mathbf{L}$; and then each inner row of $\mathbf{i}$ contains $-1$ for each subsquare entering that vertex, and $+1$ for each subsquare leaving that vertex.

In the current work's graph, the corresponding feature to an ENA node is a \emph{horizontal pocket}: two West-East paths that have the same starting vertex and ending vertex as each other, and that enclose no other West-East edges.  An example is shown in Figure~\ref{fig:level}.  There is at least one West-East path from the West cardinal vertex to the starting vertex, and at least one West-East path from the ending vertex to the East cardinal vertex.  Therefore there are two traverses that differ only in the vertices in the horizontal pocket.  Equation~(\ref{deducedlengths}) is then the difference between the corresponding traverse equations.  For the example, the relevant traverse equations are the last two in (\ref{traversesWE}).

The four linear equations in (\ref{traversesWE}) are not sufficient to fix all lengths.  For a full square dissection, the equations in (\ref{traversesNS}) also apply.  We offer no proof that the combined traverse equations are sufficient for all graphs, where \emph{sufficient} means either fixing all lengths to values consistent with all equations and with a dissection, or forming an inconsistent set of equations.  For all cases where solution has been attempted for a set of traverse equations, they have been found to be sufficient.

\subsection{Local equations}
\label{subsection:local}

Here we consider the local equations, which use coordinates $x_i$ and $y_i$ as well as lengths $L_i$.  Whenever a West-East edge connects subsquare a to subsquare b, there is a local equation: $x_a + L_a=x_b$.  Similarly, a North-South edge connecting a to e gives $y_a + L_a=y_e$.

Here we show that the local equations are equivalent to the ENA's equations, as conventionally stated in terms of $\mathbf{p}$.  We have shown that the analogy's equations can be expressed in terms of subsquare lengths, and that these equations can also be deduced from the traverse equations.  The analogy's equations have a single solution for a valid dissection, but this demonstration relies on the equations being stated in terms of $\mathbf{p}$ as in (\ref{jmatrix}).  This ensures that there are only $m-1$ unknowns in the $m-2$ equations.  This can be regarded as starting from the length-based equations and substituting (\ref{lengthpdiff}) etc.

The local equations can clearly be combined to form the traverse equations: for example, the West-East traverse through a and b can be assembled from $x_a+L_a=x_b$ and $x_b+L_b=x_E$.  (Again we note that uppercase E  and lowercase e represent cardinal and subsquare vertices respectively.)  In the ENA, $p_1$ was arbitrarily set to 0 and $p_m$ was set to either $n$ or 1.  Similarly, all subsquares connected to the West cardinal vertex have zero $x$-coordinate, and $x_E$ is $n$.  Thus we can combine these local equations to give the first traverse equation in (\ref{traversesWE}): $L_a+L_b=n$.  In practice, it was easier to produce equations in terms of normalised variables.

When the local equations from North-South edges inside a horizontal pocket are combined, it is clear that there is a single $y$-coordinate that applies to all the lower subsquares.  In the example of Figure~\ref{fig:level}, $y_f$, $y_g$ and $y_h$ are each equated to $y_c+L_c$, so they must equal each other.  This is implicit in the ENA, where edges leaving the corresponding vertex are all associated with the same $p_4$.  Finally, the expansion of lengths into $y$- or $p$-coordinates, as in (\ref{lengthpdiff}), is enforced by the North-South local equations.  Therefore, we put forward the proposition that the local equations can be rearranged to imply the equations from the vertical scan of the ENA, and similarly the horizontal scan.   It has been shown \cite{BrooksSmithStoneTutte1992} that the ENA equations, from either a horizontal or a vertical scan alone, are sufficient.  Therefore, it is conjectured that the local equations are sufficient.  More concretely, the solution procedure described in Section~\ref{equations} can detect whether a set of local equations is underdetermined, and no such sets were found in the current work.

\section{Cross-check}
\label{crosscheck}

As a check on the results, an entirely different method was used to generate all dissections for small lengths.  This simply considered all possible ways of dissecting squares of length $n$. This was tackled as an \emph{exact cover} problem: covering each of the $n^2$ unit subsquares by selecting from all possible smaller subsquares without overlap, using Knuth's Dancing Links~X algorithm \cite{Knuth1999}.

This alternative generation was used to give a complete generation of solutions of all orders, for sizes $n\leq9$.  This itself enabled a cross-check, because the results agreed with sequences in the OEIS: A045846, A221845 and A224239.  These results also agreed with the five terms of A221844 given in the OEIS, and added four new terms; details are in the Appendix.

The alternative generation was also used to give all solutions with both $n\leq18$ and $N\leq18$, and its results agreed with the main method's.  The numbers of solutions are shown in Table~\ref{crosschecktable}.

The alternative method was not intensively optimised.  The complete generation for $n=9$ required 8.5~hours on an Intel Core i7-3770 machine, as a single process.  The equivalent time for $n=8$ was 4~minutes, which suggested that optimisation would be justified before trying $n=10$.

An efficient way to extend A224239 (`Number of inequivalent ways to cut an $n \times n$ square into squares with integer sides') was to count symmetric solutions.  This was similarly regarded as an exact cover problem, selecting from symmetric collections of smaller subsquares.  These counts were combined with the published values of A045846 (`Number of distinct ways to cut an $n \times n$ square into squares with integer sides') to deduce the number of asymmetric dissections.  Each asymmetric dissection in A224239 corresponds to 8 dissections in A045846; symmetric dissections correspond to only 4, 2 or 1 dissections, depending on the symmetries.  In this method, the counts for only the symmetric solutions took less than 10~seconds for $n=10$, 1~minute for $n=11$, and 2~hours for $n=12$.  Results are given in the Appendix.

\begin{table*}[tb]
\begin{tabular}{l|r r r r r r r r r r r r}
\hline
$n$ & \multicolumn{12}{c}{Order $N$} \rule{0pt}{2.5ex} \\
\cline{2-13}
    & 7 & 8 & 9 & 10 & 11 & 12 & 13 & 14 & 15 & 16 & 17 &  \rule{0pt}{2.5ex} 18 \\
\hline
 4   & 2 & 1 &   & 4 &   &   & 3 &   &   & 1 &   &   \rule{0pt}{2.5ex} \\
 5   &   & 5 &   & 1 & 10 &   & 14 & 6 &   & 20 & 3 &   \\
 6   &   &   & 11 &   & 9 & 28 & 12 & 47 & 58 & 85 & 27 & 151   \\
 7   &   &   & 4 & 14 &   & 42 & 40 & 83 & 175 & 103 & 509 & 620 \\
 8   &   &   &   & 28 & 17 & 13 & 99 & 156 & 325 & 492 & 1291 & 1135  \\
 9   &   &   &   & 9 & 28 & 25 & 129 & 142 & 336 & 1357 & 1254 & 3256  \\
10   &   &   &   &   & 75 & 13 & 86 & 474 & 371 & 1910 & 2701 & 7150  \\
11   &   &   &   &   & 34 & 67 & 30 & 342 & 709 & 1341 & 4377 & 9712  \\
12   &   &   &   &   & 9 & 176 & 97 & 327 & 1617 & 1850 & 6622 & 15402  \\
13   &   &   &   &   & 1 & 162 & 24 & 219 & 1134 & 1759 & 7600 & 13377  \\
14   &   &   &   &   &   & 93 & 244 & 260 & 1634 & 2926 & 9271 & 21419  \\
15   &   &   &   &   &   & 27 & 432 & 244 & 781 & 4853 & 7558 & 21220  \\
16   &   &   &   &   &   & 9 & 493 & 393 & 852 & 5608 & 9943 & 32117  \\
17   &   &   &   &   &   & 2 & 216 & 309 & 594 & 4101 & 6076 & 27442  \\
18   &   &   &   &   &   &   & 187 & 1135 & 1045 & 5757 & 14362 & 36856  \\
 \hline
\end{tabular}
\caption{Number of prime dissections of order $N$ of a square of length $n$ up to symmetry.  Blanks indicate zero.  These values were found by both methods.  As well as the trivial prime dissection at $(1,1)$, single dissections exist at smaller $n$ for $(n,N)$ pairs $(2,4)$, $(3,6)$ and $(3,9)$.  An aside: those are \emph{island} non-zero values, surrounded by zeros on all four sides.  Another non-zero island is at $(6,36)$, with a zero island (surrounded by non-zeros) next to it at $(7,36)$.  These, with $(6,10)$ and $(7,11)$ seen above, are the only known islands.}
\label{crosschecktable}
\end{table*}

\section{Results}
\label{results}

The exhaustive generation was conducted for orders up to 18.  The computer time increased by a factor close to 10 for each successive order.  For example, orders 15, 16 and 17 required 0.8, 8.0 and 77~hours respectively on a four-core Intel Core i7-3770 machine, running eight processes simultaneously.  Order~18 required 18.4~days on an eight-core Intel Xeon E5-2680 machine, running sixteen processes simultaneously.

Up to symmetry, the numbers of prime dissections of order $N=1,\dots,18$ are 1, 0, 0, 1, 0, 1, 2, 6, 16, 56, 183, 657, 2277, 8813, 34178, 137578, 558734, 2285694.  This sequence (A221841) and others are listed in the Appendix.  The numbers of graphs and solutions are given in Table~\ref{orderresultstable}.

\begin{table*}[tb]
\begin{tabular}{l|r r r r}
\hline
$N$ & Graphs considered & Graphs solved & Solutions &   \rule{0pt}{2.5ex} Dissections \\
\hline
4 & 1 & 1 & 2 &  \rule{0pt}{2.5ex} 1 \\
5 & 4 & 0 & 0 & 0 \\
6 & 17 & 2 & 4 & 1 \\
7 & 89 & 5 & 17 & 2 \\
8 & 491 & 11 & 21 & 6 \\
9 & 2806 & 43 & 219 & 16 \\
10 & 16534 & 127 & 543 & 56 \\
11 & 98587 & 446 & 1711 & 183 \\
12 & 594236 & 1584 & 8637 & 657 \\
13 & 3607916 & 5761 & 32482 & 2277 \\
14 & 22046012 & 22245 & 129231 & 8813 \\
15 & 135456226 & 86138 & 571817 & 34178 \\
16 & 836535543 & 345792 & 2544446 & 137578 \\
17 & 5190532666 & 1407792 & 10677363 & 558734 \\
18 & 32348310237 & 5767035 & 47665680 & 2285694 \\
 \hline
\end{tabular}
\caption{Results for quilt order $N$.  `Graphs considered' is the number of planar imbeddings of triangulations of a square with $(N+4)$ vertices of minimum degree 4, unique up to isomorphism. `Graphs solved' is the number of these graphs that have at least one direction and colouring that leads to a solution in the sizes of a square dissection; `Solutions' is the number of these solutions.  `Dissections' is the number of prime dissections up to symmetry (A221841).  Different graphs can correspond to the same dissection, because of the alternative ways to deal with crosses.  For the same reason, a graph can have many valid ways to be directed and coloured.   For all results in the current work, all solutions from a graph correspond to the same dissection.}
\label{orderresultstable}
\end{table*}

The results confirm that the known dissections for orders up to 18 (available at, for example, \cite{Gay_at_Pegg2013}) include the largest possible.  Table~\ref{endtable} shows that, for example, the single known solution for $N=18$, $n=91$ is uniquely, up to symmetry, the largest possible at that order.  The largest size, defined as $g(N)$ in Table~\ref{endtable}, is Sequence A089047. Since dissections of order 19 are known \cite{Gay_at_Pegg2013} up to $n=120$, then the uncertain, incorrect list in \cite{CroftFalconerGuy1991} can be definitely replaced by the one in Table~\ref{ftable}.  Thus the uncertainty in Sequence A005670 has been removed -- or in fact, of course, displaced to higher $n$.

\begin{table*}[tb]
\begin{tabular}{l | c | l r r r r r r r r r r }
\hline
$N$ & $g(N)$ & \multicolumn{11}{c}{Number of dissections of order $N$, size $(g(N)-i)$}\rule{0pt}{2.5ex} \\
\cline{3-13}
   &  & $i=$ & 9 & 8 & 7 & 6 & 5 & 4 & 3 & 2 & 1 & 0\rule{0pt}{2.5ex} \\
\hline
12 &   17 & &   13 &   25 &   13 &   67 &  176 &  162 &   93 &   27 &    9 &    2 \\
13 &   23 & &  244 &  432 &  493 &  216 &  187 &  113 &   44 &   11 &    1 &    2 \\
14 &   29 & & 1231 &  798 &  696 &  561 &  288 &  152 &   53 &   35 &   17 &    4 \\
15 &   41 & &   386 &  210 &  104 &   65 &   35 &   17 &    2 &    5 &   0 &    1 \\
16 &   53 & &  231 &  127 &   81 &   62 &   27 &   16 &    6 &    3 &    1 &    1 \\
17 &   70 & &  131 &   57 &   42 &   26 &   11 &   10 &    4 &    4 &    1 &    3 \\
18 &   91 & &   71 &   35 &   33 &   16 &   14 &   12 &    1 &    3 &    3 &    1 \\
 \hline
\end{tabular}
\caption{Numbers of prime dissections of order $N$, up to symmetry, for sizes close to and equal to the largest possible, $g(N)$.  The results for $N\leq12$ can be seen in the transpose of Table~\ref{crosschecktable}.
}
\label{endtable}
\end{table*}

\begin{table*}[tb]
\begin{tabular}{r  c c c c c c c c c c c c}
\hline
$f(n)=$ & 1 & 4 & 6 & 7 & 8 & 9     & 10   & 11 & 12 & 13 & 14\rule{0pt}{2.5ex} \\
$n=$     & 1 & 2 & 3 & 4 & 5 & 6,7 &  8,9 & 10--13 & 14--17 & 18--23 & 24--29\rule{0pt}{2.5ex} \\
\hline
\end{tabular}
\begin{tabular}{r  c c c c c c }
$f(n)=$ & 15 & 16 & 17 & 18 & 19\rule{0pt}{2.5ex} \\
$n=$     & 30--39,41 & 40,42--53 & 54--70 & 71--91 & 92--120,122,126,\ldots?\rule{0pt}{2.5ex} \\
\hline
\end{tabular}
\caption{The smallest possible order, $f(n)$, of a prime dissection of a square of size $n$.}
\label{ftable}
\end{table*}

\appendix

\section{Confirmed and extended sequences in OEIS}
\label{appendix}

This appendix lists values from several sequences in the OEIS \cite{OEIS} with previously known values that have been confirmed in the current work (or, conversely, have been used to check the current work).  \textit{Also, new (or newly definite) values are shown in italics}.

These sequences, as quoted in \cite{OEIS}, include the trivial dissection of a square into itself (except where `smaller' is specified, in A018835 and A211302).

Sequences A045846 and A221845 have known values beyond those listed here.  In particular, the values of A045846$(n)$ for $n \in \{10, 11, 12\}$ (1500957422222, 790347882174804 and 781621363452405930) were used to confirm and extend A224239, as described in Section~\ref{crosscheck}.

{\footnotesize{ \raggedright{
\begin{itemize}
\item A005670 (Mrs. Perkins's quilt: smallest coprime dissection of $n \times n$ square) for $n=1,\dots,120$:
    1,  4,  6,  7,  8,  9,  9, 10, 10, 11, 
   11, 11, 11, 12, 12, 12, 12, 13, 13, 13, 
   13, 13, 13, 14, 14, 14, 14, 14, 14, 15, 
   15, 15, 15, 15, 15, 15, 15, 15, 15, \textit{16}, 
   15, \textit{16, 16, 16, 16, 16, 16, 16, 16, 16, 
   16, 16, 16, 17, 17, 17, 17, 17, 17, 17, 
   17, 17, 17, 17, 17, 17, 17, 17, 17, 17, 
   18, 18, 18, 18, 18, 18, 18, 18, 18, 18, 
   18, 18, 18, 18, 18, 18, 18, 18, 18, 18, 
   18, 19, 19, 19, 19, 19, 19, 19, 19, 19, 
   19, 19, 19, 19, 19, 19, 19, 19, 19, 19, 
   19, 19, 19, 19, 19, 19, 19, 19, 19, 19}.
\item A045846 (Number of distinct ways to cut an $n \times n$ square into squares with integer sides) for $n=1,\dots,9$:
   1, 2, 6, 40, 472, 10668, 450924, 35863972, 5353011036.
\item A089046 (Least edge-length of a square dissectable into at least $N$ squares in the Mrs Perkins's quilt problem) for $N=1,\dots,19$:
    1, 2, 2, 2, 3, 3, 4, 5, 6, 8, 10, 14, 18, 24, 30, \textit{40, 54, 71, 92}.
\item A089047 (Greatest edge-length of a square dissectable into up to $N$ squares in Mrs Perkins's quilt problem) for $N=1,\dots,18$:
   1, 1, 1, 2, 2, 3, 4, 5, 7, 9, 13, 17, 23, 29, \textit{41, 53, 70, 91}.
\item A018835 (Minimal number of smaller integer-sided squares that tile an $n \times n$ square) for $n=2,\dots,126$:
   4, 6, 4, 8, 4, 9, 4, 6, 4, 11, 4, 11, 4, 6, 4, 12, 4, 13, 4, 6, 4, 13, 4, 8, 4, 6, 4, 14, 4, 15, 4, 6, 4, 8, 4, 15, 4, 6, 4, 15, 4, 16, 4, 6, 4, 16, 4, 9, 4, 6, 4, 16, 4, 8, 4, 6, 4, 17, 4, 17, 4, 6, 4, 8, 4, \textit{17, 4, 6, 4, 18, 4, 18, 4, 6, 4, 9, 4, 18, 4, 6, 4, 18, 4, 8, 4, 6, 4, 18, 4, 9, 4, 6, 4, 8, 4, 19, 4, 6, 4, 19, 4, 19, 4, 6, 4, 19, 4, 19, 4, 6, 4, 19, 4, 8, 4, 6, 4, 9, 4, 11, 4, 6, 4, 8, 4}.
\item A211302 (Minimal number of smaller integer-sided squares that tile a $p \times p$ square, where $p=i$-th prime) for $i=1,\dots,30$: 
   4, 6, 8, 9, 11, 11, 12, 13, 13, 14, 15, 15, 15, 16, 16, 16, 17, 17, \textit{17, 18, 18, 18, 18, 18, 19, 19, 19, 19, 19, 19}.
\item A221841 (Number of ways to dissect a square into $N$ squares up to symmetry) for $N=1,\dots,18$: 
   1, 0, 0, 1, 0, 1, 2, 6, 16, \textit{56, 183, 657, 2277, 8813, 34178, 137578, 558734, 2285694}.
\item A221842 (Number of ways to dissect a square into $N$ squares) for $N=1,\dots,18$: 
   1, 0, 0, 1, 0, 4, 8, 36, 105, \textit{384, 1340, 4975, 17676, 69052, 270716, 1093218, 4455047, 18246018}.
\item A221844 (Number of prime dissections of an $n \times n$ square into integer-sided squares up to symmetry) for $n=1,\dots,9$:
   1, 1, 2, 11, 76, \textit{1490, 56977, 4495010, 669203525}.
\item A221845 (Number of prime dissections of an $n \times n$ square into integer-sided squares) for $n=1,\dots,9$:
   1, 1, 5, 38, 471, 10661, 450923, 35863932, 5353011030.
\item A224239 (Number of inequivalent ways to cut an $n \times n$ square into squares with integer sides) for $n=1,\dots,12$:
   1, 2, 3, 13, 77, 1494, 56978, 4495023, 669203528, 187623057932, \textit{98793520541768, 97702673827558670}.
\item A226978 (Number of ways to cut an $n \times n$ square into squares with integer sides, reduced for symmetry, where the orbits under the symmetry group of the square, D4, have 1 element) for $n=1,\dots,12$:
   1, 2, 2, 4, 4, 12, 8, \textit{44, 32, 228, 148, 1632}.
\item A226979 (Number of ways to cut an $n \times n$ square into squares with integer sides, reduced for symmetry, where the orbits under the symmetry group of the square, D4, have 2 elements) for $n=1,\dots,12$:
   0, 0, 0, 2, 2, 24, 36, \textit{344, 504, 7657, 11978, 289829}.
\item A226980 (Number of ways to cut an $n \times n$ square into squares with integer sides, reduced for symmetry, where the orbits under the symmetry group of the square, D4, have 4 elements) for $n=1,\dots,12$:
   0, 0, 1, 6, 26, 264, 1157, \textit{23460, 153485, 6748424, 70521609, 6791578258}.
\item A226981 (Number of ways to cut an $n \times n$ square into squares with integer sides, reduced for symmetry, where the orbits under the symmetry group of the square, D4, have 8 elements) for $n=1,\dots,12$:
   0, 0, 0, 1, 45, 1194, 55777, \textit{4471175, 669049507, 187616301623, 98793450008033, 97702667035688951}.
\end{itemize}

The following new sequences have been added to the OEIS as results of the current work:
\begin{itemize}
\item A232484 (Number of size collections in prime `Mrs. Perkins's Quilt' dissections of integer-sided squares into $N$ squares) for $N=1,\dots,18$:
   \textit{1, 0, 0, 1, 0, 1, 1, 2, 4, 7, 18, 40, 119, 323, 1100, 3594, 13068, 47444}.
\item A240120 (Number of inequivalent ways to cut an $n \times n$ square into squares with integer sides, such that the dissection has reflective symmetry in both diagonals and no other reflective symmetries) for $n=1,\dots,12$:
   \textit{0, 0, 0, 1, 1, 9, 19, 121, 275, 2489, 7217, 86775}.
\item A240121 (Number of inequivalent ways to cut an $n \times n$ square into squares with integer sides, such that the dissection has two reflective symmetries in axes parallel to the sides, and no other reflective symmetries) for $n=1,\dots,12$:
   \textit{0, 0, 0, 1, 0, 13, 5, 183, 75, 4408, 1501, 180324}.
\item A240122 (Number of inequivalent ways to cut an $n \times n$ square into squares with integer sides, such that the dissection has 90-degree rotational symmetry and no reflective symmetry) for $n=1,\dots,12$:
   \textit{0, 0, 0, 0, 1, 2, 12, 40, 154, 760, 3260, 22730}.
\item A240123 (Number of inequivalent ways to cut an $n \times n$ square into squares with integer sides, such that the dissection has a reflective symmetry in one diagonal, but no other symmetries) for $n=1,\dots,12$:
   \textit{0, 0, 1, 3, 19, 107, 847, 8647, 119835, 2255123, 58125783, 2050662011}.
\item A240124 (Number of inequivalent ways to cut an $n \times n$ square into squares with integer sides, such that the dissection has 180-degree rotational symmetry, but no other symmetries) for $n=1,\dots,12$:
   \textit{0, 0, 0, 0, 2, 19, 109, 1781, 13660, 397689, 5368943, 289864745}.
\item A240125 (Number of inequivalent ways to cut an  $n \times n$ square into squares with integer sides, such that the dissection has one reflective symmetry in an axis parallel to a side, but no other symmetries) for $n=1,\dots,12$:
   \textit{0, 0, 0, 3, 5, 138, 201, 13032, 19990, 4095612, 7026883, 4451051502}.
\end{itemize}
} } }

\medskip
\noindent{\bf Acknowledgement:}  I am grateful for the helpful comments from the reviewers at \emph{Discrete Mathematics}, and also from Dr.~Matthew Fayers at Queen Mary, University of London, UK.

\end{document}